\documentclass[a4paper,12pt]{article}
\usepackage{a4wide}
\usepackage{theorem}
\usepackage{amsmath}
\usepackage{array}
\usepackage{amssymb}
\usepackage{amsfonts}
\usepackage[english]{babel}
\usepackage{epsf}
\usepackage{epsfig}
\usepackage{graphicx}

\newtheorem{theo}{\indent Theorem \newline}[section]
\newtheorem{defi}[theo]{\indent Definition\newline}
\newtheorem{rem}[theo]{\noindent Remark}
\newtheorem{prop}[theo]{\indent Proposition\newline}
\newtheorem{lemma}[theo]{\indent Lemma\newline}
\newtheorem{cor}[theo]{\indent Corollary \newline}

 \def\N{{\mathbb{N}}}
\def\Z{{\mathbb{Z}}}

\def\R{{\mathbb{R}}}
\def\C{{\mathbb{C}}}
\def\H{{\mathbb{H}}}

\usepackage{eufrak}

\newlength{\indentation}%
\makeatletter
\newcommand\@makefntextsans[1]{%
    \parindent 0em%
    \noindent%
    \hb@xt@0em{\hss}%
    #1}
\def\footnotetextsans{%
     \@ifnextchar [\@xfootnotenextsans%
       {\@footnotetextsans}}
\def\@xfootnotenextsans[#1]{%
  \begingroup%
     \csname c@\@mpfn\endcsname #1\relax%
  \endgroup%
  \@footnotetextsans}
\long\def\@footnotetextsans#1{\insert\footins{%
    \reset@font\footnotesize%
    \interlinepenalty\interfootnotelinepenalty%
    \splittopskip\footnotesep%
    \splitmaxdepth \dp\strutbox \floatingpenalty \@MM%
    \hsize\columnwidth \@parboxrestore%
    \color@begingroup%
      \@makefntextsans{%
        \rule\z@\footnotesep\ignorespaces#1\@finalstrut\strutbox}
    \color@endgroup}}
\makeatother

\begin{document}

\title{Open Gromov-Witten invariants in dimension four}
\author{Jean-Yves Welschinger}
\maketitle

\makeatletter\renewcommand{\@makefnmark}{}\makeatother
\footnotetextsans{Keywords: holomorphic discs, Gromov-Witten invariants.}
\footnotetextsans{AMS Classification : 53D45.
}

{\bf Abstract:}

Given a closed orientable Lagrangian surface $L$ in a closed symplectic four-manifold $(X , \omega)$
together with a relative homology class $d \in H_2 (X , L ; \Z)$ with vanishing boundary in $H_1 (L ; \Z)$, we prove that
the algebraic number of $J$-holomorphic discs with boundary on $L$, homologous to $d$ and passing
through the adequate number of points neither depends on the choice of the points nor on the generic choice
of the almost-complex structure $J$. We furthermore get analogous open Gromov-Witten invariants by counting,
for every non-negative integer $k$, unions of $k$ discs instead of single discs. 

\section*{Introduction}
 
 Let $(X , \omega)$ be a closed connected symplectic four-manifold. Let $L \subset X$ be a closed Lagrangian surface
 which we mainly assume to be orientable. We denote by $\mu_L \in H^2 (X, L ; \Z)$ its Maslov class, that is the obstruction to 
 extend $TL$ as a Lagrangian subbundle of $TX$. We denote by ${\cal J}_\omega$ the space of almost-complex
 structures of class $C^l$ tamed by $\omega$, where $l \gg 1$ is a fixed integer. Let $d \in H_2 (X, L ; \Z)$ be
 such that $\mu_L (d) > 0$ and $r,s \in \N$ such that $r + 2s = \mu_L (d) - 1$. Let $\underline{x} \subset L^r$
 (resp. $\underline{y} \subset (X \setminus L)^s$) be a collection of $r$ (resp. $s$) distinct points. Then, for every generic
 choice of $J \in {\cal J}_\omega$, $X$ contains only finitely many $J$-holomorphic discs with boundary on $L$, homologous to $d$ and
 which pass through $\underline{x} \cup \underline{y}$. These discs are all immersed and we denote by
 ${\cal M}_d (\underline{x}, \underline{y} , J)$ their finite set. For every $D \in {\cal M}_d (\underline{x}, \underline{y} ; J)$,
 we denote by $m(D) = [\stackrel{\circ}{D}] \circ [L] \in \Z/2\Z$ the intersection index between the interior of $D$ and the
 surface $L$. We then set
 $$GW_d^r(X,L ; \underline{x}, \underline{y} , J) = \sum_{D \in {\cal M}_d (\underline{x}, \underline{y} ; J)} (-1)^{m(D)} \in \Z.$$
 Our main result is the following (see Theorem \ref{theo1disc}). 
 \begin{theo}
 \label{theo0}
 Assume that $L$ is orientable and that $\partial d = 0 \in H_1 (L ; \Z)$. Then, the integer 
 $GW_d^r(X,L ; \underline{x}, \underline{y} , J)$ neither depends on the choice of $\underline{x}, \underline{y}$ nor on the generic choice of $J$. $\square$
 \end{theo}
 This result thus provides an integer valued invariant which we can denote without ambiguity
 by $GW_d^r(X,L )$, providing a relative analog to the genus zero Gromov-Witten invariants in dimension four. 
 Note that some open Gromov-Witten invariants have already been defined by C.C. Liu and M. Katz in the presence of an action of the circle, see \cite{KatzLiu}, \cite{Liu},
 and by myself when $L$ is fixed by an antisymplectic involution, see \cite{WelsCRAS}, \cite{WelsInvent}, \cite{WelsICM}
 or also \cite{Cho}, \cite{Solomon}. 
 When $X$ is six-dimensional, some open Gromov-Witten invariants have been defined by K. Fukaya \cite{Fuk} and V. Iacovino \cite{Iac} in the case
 of Calabi-Yau six-manifolds and by myself \cite{Welsopen6} in the absence of Maslov zero discs. 
 
 When $L$ is a Lagrangian sphere
 fixed by such an antisymplectic involution, the invariant $\chi^d_r$ introduced in \cite{WelsCRAS1}, \cite{WelsInvent} is computed in terms of
 $GW_d^r(X,L )$, see Lemma \ref{lemmadivide}. We deduce as a consequence that $2^{s-1}$ divides $\chi^d_r$, improving 
 a congruence already obtained in \cite{WelsSFT} using symplectic field theory. 
 When the genus of $L$ is greater than one, the invariant $GW_d^r(X,L )$ vanishes, see Proposition \ref{propSFT}, since one can find a generic
 almost complex structure $J$ for which the set ${\cal M}_d (\underline{x}, \underline{y} ; J)$ is empty. 
 I could not get any result when $L$ is not orientable, except a result modulo two when $L$ is homeomorphic to a real projective plane, see Theorem \ref{theo1discno}. 
 
We also obtain an analog of Theorem \ref{theo0} by counting 
unions of $k$ discs, $k > 0$, instead of single discs. More precisely, let $k >0$ be such that $\mu_L (d) \geq k$ and assume now
that $r + 2s = \mu_L (d) - k$. For every generic
 choice of $J \in {\cal J}_\omega$, $X$ only contains finitely many unions of $k$ $J$-holomorphic discs with boundary on $L$, total homology class $d$ and
 which pass through $\underline{x} \cup \underline{y}$. These discs are all immersed and we denote by
 ${\cal M}_{d,k} (\underline{x}, \underline{y} , J)$ their finite set. For every ${ D} =D_1 \cup \dots \cup D_k \in {\cal M}_{d,k} (\underline{x}, \underline{y} ; J)$,
 we denote by $m(D) = \sum_{i=1}^{k} m(D_i) \in \Z/2\Z$ and set
 $$GW_{d,k}^r(X,L ; \underline{x}, \underline{y} , J) = \sum_{{ D} \in {\cal M}_{d,k} (\underline{x}, \underline{y} ; J)} (-1)^{m({ D})} \in \Z.$$
 We then get (see Theorem \ref{theokdiscs}). 
 \begin{theo}
 \label{theo0'}
 Assume that $L$ is orientable and that $\partial d = 0 \in H_1 (L ; \Z)$. Then, the integer $GW_{d,k}^r(X,L ; \underline{x}, \underline{y} , J)$
 neither depends on the choice of $\underline{x}, \underline{y}$ nor on the generic choice of $J$. $\square$
 \end{theo}
Note that the individual discs involved in the definition of this k-discs open Gromov-Witten invariant $GW_{d,k}^r(X,L) \in \Z$ are no more subject
to have trivial boundary in homology. 

The first part of this paper is devoted to notations and generalities on moduli spaces of pseudo-holomorphic discs in any dimensions. 
We introduce the numbers $GW_d^r(X,L )$, $GW_{d,k}^r(X,L)$ and prove their invariance in the second paragraph. \\

{\bf Acknowledgements:}

This paper is mostly based on an idea I got shortly after my stay at the ETH Z\"urich in the Fall $2010$. I wish to acknowledge both
ETH Z\"urich for its hospitality and Paul Biran for the many fruitful discussions we had. The research leading to these results has received funding
from the European Community's Seventh Framework Progamme ([FP7/2007-2013] [FP7/2007-2011]) under
grant agreement $\text{n}\textsuperscript{o}$ [258204], as well as from the French Agence nationale de la recherche, ANR-08-BLAN-0291-02. 
Finally, I am grateful to Erwan Brugall\'e for pointing out the reference \cite{AbraBert} to me, which is related to his work \cite{BruPui} with 
N. Puignau. It led me to revise a first version of this paper.

\section{Pseudo-holomorphic discs with boundary on a Lagrangian submanifold}
\subsection{Moduli spaces of simple discs}

Let $\Delta = \{ z \in \C \; \vert \; \vert z \vert \leq 1 \}$ be the closed complex unit disc. We denote by
${\cal P} (X,L) = \{ (u,J) \in C^{1} (\Delta , X) \times {\cal J}_\omega \; \vert \; u(\partial \Delta) \subset L \text{ and } du + J\vert_u \circ du \circ j_{st} = 0\}$
the space of pseudo-holomorphic maps from $\Delta$ to the pair $(X,L)$, where $ j_{st}$ denotes the standard complex structure of $\Delta$.
Note that $J$ being of class $C^l$, the regularity of such pseudo-holomorphic maps $u$ is actually more than $C^l$, see \cite{McDSal}. 
More generally, for every $r,s \in \N$, we denote by ${\cal P}_{r,s} (X,L) = \{ ((u,J), \underline{z}, \underline{\zeta}) \in {\cal P} (X,L) \times ((\partial \Delta)^r \setminus
\text{diag}_{\partial \Delta}) \times ((\stackrel{\circ}{\Delta})^s \setminus
\text{diag}_{\Delta}) \}$, where $\text{diag}_{\partial \Delta} = \{ (z_1, \dots , z_r) \in (\partial \Delta)^r \; \vert \; \exists i \neq j , z_i = z_j \}$ and
$\text{diag}_{ \Delta} = \{ (\zeta_1, \dots , \zeta_s) \in \Delta^s \; \vert \; \exists i \neq j , \zeta_i = \zeta_j \}$.

Following \cite{Laz}, \cite{KwonOh}, \cite{BiranCor}, we define
\begin{defi}
\label{defsimple}
A pseudo-holomorphic map $u$ is said to be simple iff there is a dense open subset $\Delta_{inj} \subset \Delta$ such that
$\forall z \in \Delta_{inj}, u^{-1} (u(z)) = \{ z \}$ and $du\vert_z \neq 0$.
\end{defi}
Recall for instance that the map $z \in \C \mapsto z^3 \in \C$ restricted to the upper half plane $\H \subset \C$ induces,
after composition by a biholomorphism $\Delta \to \overline{\H} \subset \C P^1$, a holomorphic map $u : \Delta \to (\C P^1 , \R P^1)$
which is not simple in the sense of Definition \ref{defsimple}, though it is somewhere injective, compare \cite{McDSal}. 

We denote by ${\cal P}_{r,s}^* (X,L) $ the subset of simple elements of ${\cal P}_{r,s} (X,L) $. It is a separable Banach manifold which is
naturally embedded as a submanifold of class $C^{l-k}$ of the space $W^{k,p} (\Delta , X)  \times {\cal J}_\omega $ for every
$1 \ll k \ll l$ and $p > 2$, see Proposition $3.2$ of \cite{McDSal}. 

For every $d \in H_2 (X , L ; \Z)$, we denote by ${\cal P}_{r,s}^d (X,L) = \{ (u,J) \in {\cal P}_{r,s}^* (X,L)  \; \vert \; u_* [\Delta ] = d \}$ and by
${\cal M}_{r,s}^d (X,L)  = {\cal P}_{r,s}^d (X,L) /\text{Aut}(\Delta) $, where $\text{Aut}(\Delta)$ is the group of biholomorphisms of $\Delta$
which acts by composition on the right. The latter is equipped with a projection $\pi : [u,J, \underline{z}, \underline{\zeta}] \in {\cal M}_{r,s}^d (X,L) \mapsto
J \in {\cal J}_\omega$ and an evaluation map $\text{eval} :  [u,J, \underline{z}, \underline{\zeta}] \in {\cal M}_{r,s}^d (X,L) \mapsto (u(\underline{z}), u(\underline{\zeta}))
\in L^r \times X^s$. 

We recall the following classical result due to Gromov (see \cite{Gromov}, \cite{McDSal}, \cite{FOOO}).
\begin{theo}
\label{theodim}
For every closed Lagrangian submanifold $L$ of a $2n$-dimensional closed symplectic manifold $(X , \omega)$
and for every $d \in H_2 (X , L ; \Z)$, $r,s \in \N$, the space ${\cal M}_{r,s}^d (X,L)$ is a separable Banach manifold and the projection
$\pi :  {\cal M}_{r,s}^d (X,L) \to {\cal J}_\omega$ is Fredholm of index $\mu_L (d) + n - 3 + r +2s$. $\square$
\end{theo}
Note that from Sard-Smale's theorem \cite{Smale}, the set of regular values of $\pi$ is dense of the second category. As a consequence,
for a generic choice of $J \in  {\cal J}_\omega$, the moduli space ${\cal M}_{0,0}^d (X,L ; J) = \pi^{-1} (J)$ is a manifold of dimension
$\mu_L (d) + n - 3 $ as soon as it is not empty. Likewise, if $\underline{x}$ (resp.  $\underline{y}$) is a set of $r$ (resp. $s$) distinct points
of $L$ (resp. $X \setminus L$), then ${\cal M}_{r,s}^d (X,L ; \underline{x} , \underline{y}, J) =( \pi \times \text{eval})^{-1} (J, \underline{x} , \underline{y})$
is a manifold of dimension $\mu_L (d) + n - 3 -(n-1)(r+2s)$. We denote by ${\cal M}_{r,s}^d (X,L ;  \underline{x} , \underline{y})$ the preimage
$ \text{eval}^{-1} ( \underline{x} , \underline{y})$ and with a slight abuse by $\pi : [u,J, \underline{z}, \underline{\zeta}] \in {\cal M}_{r,s}^d (X,L;  \underline{x} , \underline{y}) \mapsto
J \in {\cal J}_\omega$ the Fredholm projection of index $\mu_L (d) + n - 3 -(n-1)(r+2s)$.

Recall also that the tangent bundle to the space ${\cal P}_{r,s}^* (X,L) $ writes, for every $ (u,J) \in {\cal P}_{r,s}^* (X,L)$,
$$T_{(u,J)} {\cal P}_{r,s}^* (X,L) = \{ (v , \stackrel{.}{J}) \in \Gamma^1 (\Delta, u^*TX) \times T_J {\cal J}_\omega \; \vert \; Dv +  \stackrel{.}{J} \circ du \circ j_{st} = 0
\text{ and } v\vert_{\partial \Delta} \subset u^* TL \},$$
where $D$ is the Gromov operator defined for every $v \in \Gamma^1 (\Delta, u^*TX)$ by the formula 
$Dv = \nabla v + J \circ  \nabla v \circ  j_{st} + \nabla_v J \circ du \circ  j_{st} \in \Gamma^0 (D, \Lambda^{0,1} D \otimes u^* TX)$ for any torsion free
connection $\nabla$ on $TX$. This operator
induces an operator $\overline{D}$ on the normal sheaf ${\cal N}_u = u^* TX / du(T\Delta)$, see formula $1.5.1$ of \cite{IvashShev},
and we denote by $H_D^0 (\Delta , {\cal N}_u ) \subset \Gamma^1 (\Delta, u^*TX, u^*TL)/du(\Gamma^1 (\Delta, T\Delta , T\partial \Delta))$
the kernel of this operator $\overline{D}$, see \cite{IvashShev} or \S $1.4$ of \cite{WelsInvent}. Likewise, we denote by
$H_D^0 (\Delta , {\cal N}_{u, - \underline{z}, -\underline{\zeta}}) $ the kernel of $\overline{D}$ restricted to elements which vanish at the points $\underline{z}$ and $ \underline{\zeta}$.
\begin{prop}
\label{proptangent}
Under the hypothesis of Theorem \ref{theodim}, let $\underline{x}$ (resp.  $\underline{y}$) be a set of $r$ (resp. $s$) distinct points
of $L$ (resp. $X \setminus L$) and let $J $ be a generic element of $ {\cal J}_\omega$. Then, at every point $[u,J, \underline{z}, \underline{\zeta}] $
of $ {\cal M}_{r,s}^d (X,L; J, \underline{x} , \underline{y}) $ the tangent space $T_{[u,J, \underline{z}, \underline{\zeta}] } {\cal M}_{r,s}^d (X,L; J, \underline{x} , \underline{y}) $
is isomorphic to $H_D^0 (\Delta , {\cal N}_{u, - \underline{z}, -\underline{\zeta}}) $. $\square$
\end{prop}
This classical result is proved in \cite{McDSal} or  \cite{IvashShev} for example (compare \S $1.8$ of \cite{WelsInvent}). 

Recall finally that the moduli space ${\cal M}_{r,s}^d (X,L; J, \underline{x} , \underline{y}) $ given by Proposition \ref{proptangent} is not in general compact
for two reasons. Firstly, a sequence of elements of ${\cal M}_{r,s}^d (X,L; J, \underline{x} , \underline{y}) $  may converge to a pseudo-holomorphic
disc which is not simple in the sense of Definition \ref{defsimple}, see \S \ref{subsecsimple}. Secondly, a sequence of elements of ${\cal M}_{r,s}^d (X,L; J, \underline{x} , \underline{y}) $  
may converge to a pseudo-holomorphic curve which no more admits a parameterization by a single disc $\Delta$, in the same way as a sequence of smooth
plane conics may converge to a pair of distinct lines. However, the latter phenomenon is well understood by the following Gromov compactness' theorem.
\begin{theo}
\label{theocompact}
Under the hypothesis of Proposition \ref{proptangent}, every sequence of elements of $ {\cal M}_{r,s}^d (X,L; J, \underline{x} , \underline{y}) $ has a subsequence
which converges in the sense of Gromov to a stable $J$-holomorphic disc. $\square$
\end{theo}
A proof of Theorem \ref{theocompact} as well as the definitions of stable discs and convergence in the sense of Gromov can
be found in \cite{Frauen}. 

\subsection{Theorem of decomposition into simple discs}
\label{subsecsimple}

We recall in this paragraph the theorem of decomposition into simple discs established by Kwon-Oh and Lazzarini, see  \cite{KwonOh}, \cite{Laz}, \cite{Laz1}. 
\begin{theo}
\label{theodecomp}
Let $L$ be a closed Lagrangian submanifold of a $2n$-dimensional closed symplectic manifold $(X , \omega)$. Let $u : (\Delta , \partial \Delta) \to (X,L)$
be a non-constant pseudo-holomorphic map. Then, there exists a graph ${\cal G} (u)$ embedded in $\Delta$ such that $\Delta \setminus {\cal G} (u)$
has only finitely many connected components. Moreover, for every connected component ${\cal D} \subset \Delta \setminus {\cal G} (u)$, there exists a
surjective map $\pi_{\overline{\cal D}} : \overline{\cal D} \to \Delta$, holomorphic on ${\cal D}$ and continuous on $\overline{\cal D}$, as well as a simple
pseudo-holomorphic map $u_{\cal D} : \Delta \to X$ such that $u\vert_{\overline{\cal D}} = u_{\cal D} \circ \pi_{\overline{\cal D}}$. The map $\pi_{\overline{\cal D}}$
has a well defined degree $m_{\cal D} \in \N$, so that $u_* [\Delta] = \sum_{\cal D} m_{\cal D} (u_{\cal D})_* [\Delta] \in H_2 (X,L ; \Z)$, the sum being taken
over all connected components ${\cal D}$ of $ \Delta \setminus {\cal G} (u)$. $\square$
\end{theo}
The graph $ {\cal G} (u)$ given by Theorem \ref{theodecomp} is called the frame or non-injectivity graph, see \cite{Laz}, \cite{Laz1} (or \S $3.2$ of \cite{BiranCor}) for its definition. 

\subsection{Pseudo-holomorphic discs in dimension four}
\label{subsecdimensionfour}

We assume in this paragraph that the ambient closed symplectic manifold  $(X , \omega)$ is of dimension four and recall several facts
specific to this dimension. 

\begin{prop}
\label{propimmersion}
Let $L$ be a closed Lagrangian surface of a closed connected symplectic four-manifold $(X , \omega)$. Let $d \in H_2 (X , L ; \Z)$ be such that $\mu_L (d) >0$ and
 $r,s \in \N$ such that $r + 2s = \mu_L (d) - 1$. Let $\underline{x}$ (resp.  $\underline{y}$) be a set of $r$ (resp. $s$) distinct points
of $L$ (resp. $X \setminus L$). Then, the critical points of the projection $\pi : [u,J, \underline{z}, \underline{\zeta}] \in {\cal M}_{r,s}^d (X,L;  \underline{x} , \underline{y}) \mapsto
J \in {\cal J}_\omega$ are those for which $u$ is not an immersion. 
\end{prop}

{\bf Proof:}

This Proposition \ref{propimmersion} is analogous to Lemma $2.13$ of \cite{WelsInvent} and follows from the automatic transversality in
dimension four, see Theorem $2$ of \cite{HLS}. From Theorem \ref{theodim}, $\pi$ is of vanishing index whereas from Proposition
\ref{proptangent}, its kernel is isomorphic to $H_D^0 (\Delta , {\cal N}_{u, - \underline{z}, -\underline{\zeta}}) $. When $u$ is not an immersion,
the sheaf ${\cal N}_{u, - \underline{z}, -\underline{\zeta}} $ contains a skyscraper part carried by its critical points and which contributes to 
the kernel $H_D^0 (\Delta , {\cal N}_{u, - \underline{z}, -\underline{\zeta}}) $, so that $[u,J, \underline{z}, \underline{\zeta}] $ is indeed a critical point of $\pi$.
When $u$ is an immersion, this normal sheaf is the sheaf of sections of a bundle of Maslov index $-1$. From Theorem $2$ of \cite{HLS},
$H_D^0 (\Delta , {\cal N}_{u, - \underline{z}, -\underline{\zeta}}) $ is then reduced to $\{ 0 \}$. $\square$

\begin{prop}
\label{propnondeg}
Under the hypothesis of Proposition \ref{propimmersion} :

1) If $J \in {\cal J}_\omega$ is generic, then all elements of ${\cal M}_{r,s}^d (X,L; J,  \underline{x} , \underline{y})$ are immersed discs. 

2) If $t \in [0,1] \mapsto J_t \in {\cal J}_\omega$ is a generic path, then every element of $\cup_{t \in [0,1]} {\cal M}_{r,s}^d (X,L; J_t ,  \underline{x} , \underline{y})$
which is not immersed has a unique ordinary cusp which is on $\partial \Delta$. The latter are non-degenerated critical points of $\pi$. 
\end{prop}

The ordinary cusps given by Proposition \ref{propnondeg} are by definition modeled on the map $t \in \R \mapsto (t^2 , t^3) \in \R^2$ at the
neighborhood of the origin. \\

{\bf Proof:}

The first part follows from Proposition \ref{propimmersion} and Sard-Smale's theorem \cite{Smale}. 
The second part is proven exactly in the same way as the first part of Proposition $2.7$ of \cite{WelsInvent}
and the fact that these critical points are non-degenerated follows along the same lines as Lemma $2.13$
of \cite{WelsInvent}. We do not reproduce these proofs here. $\square$\\

\section{Open Gromov-Witten invariants in dimension four}
\subsection{One-disc open Gromov-Witten invariants}
\label{subsec1disc}

Let $(X , \omega)$ be a closed connected symplectic four-manifold.
Let $L \subset X$ be a closed Lagrangian surface, of Maslov class $\mu_L \in H^2 (X,L ; \Z)$. 
Let $d \in H_2 (X , L ; \Z)$ be such that $\mu_L (d) >0$ and
 $r,s \in \N$ such that $r + 2s = \mu_L (d) - 1$. Let $\underline{x}$ (resp.  $\underline{y}$) be a set of $r$ (resp. $s$) distinct points
of $L$ (resp. $X \setminus L$). For every $J \in {\cal J}_\omega$ generic, the moduli space ${\cal M}_{r,s}^d (X,L; J,  \underline{x} , \underline{y})$
is then finite and consists only of immersed discs, see Proposition \ref{propnondeg}. We denote, for every 
$ [u,J, \underline{z}, \underline{\zeta}] \in {\cal M}_{r,s}^d (X,L; J, \underline{x} , \underline{y})$, by $m(u)$ the intersection index
$[u( \stackrel{\circ}{\Delta})] \circ [L] \in \Z/2\Z$, where $\stackrel{\circ}{\Delta} = \{ z \in \C \; \vert \; \vert z \vert < 1 \}$ denotes the interior of the
disc $\Delta$. In fact, $J$ being generic, the intersection $u(\stackrel{\circ}{\Delta}) \cap L$ is transversal, so that the intersection index $m(u)$
coincides with the parity of the number of intersection points between $u(\stackrel{\circ}{\Delta} )$ and $L$, compare \S $2.1$ of \cite{WelsInvent}. 

We then set 
$$GW_d^r(X,L ; \underline{x}, \underline{y} , J) = \sum_{ [u,J, \underline{z}, \underline{\zeta}] \in {\cal M}_{r,s}^d (X,L; J,  \underline{x} , \underline{y})} (-1)^{m(u)} \in \Z.$$
 \begin{theo}
 \label{theo1disc}
 Let $(X , \omega)$ be a closed connected symplectic four-manifold and $L \subset X$ be a closed Lagrangian surface
 which we assume to be orientable. Let $d \in H_2 (X , L ; \Z)$  be such that $\mu_L (d) >0$
 and $\partial d = 0 \in H_1 (L ; \Z)$. Let $r,s \in \N$ be such that $r + 2s = \mu_L (d) - 1$ and 
 $\underline{x}$ (resp.  $\underline{y}$) be a collection of $r$ (resp. $s$) distinct points
of $L$ (resp. $X \setminus L$). Let finally $J \in {\cal J}_\omega$ be generic.
 Then, the integer $GW_d^r(X,L ; \underline{x}, \underline{y} , J)$
 neither depends on the choice of $\underline{x}, \underline{y}$ nor on the generic choice of $J$. 
 \end{theo}
 
 Before proving Theorem \ref{theo1disc} in section \ref{subsecproof1disc}, let us first formulate further results or consequences.
 The invariant provided by Theorem \ref{theo1disc} can be denoted without ambiguity by $GW_d^r(X,L) \in  \Z$. 
 As the usual Gromov-Witten invariants, it also does not change under deformation of the symplectic form $\omega$,
 whereas it indeed depends in general on $d$ and $r$, compare \S $3$ of \cite{WelsInvent}. Note that Theorem \ref{theo1disc} 
 also holds true when $(X , \omega)$  is convex at infinity. Moreover, if $X$ does not contain any Maslov zero pseudo-holomorphic
 disc with boundary on $L$, for instance if $L$ is monotone, then assuming in Theorem \ref{theo1disc} that $\partial d = 0 \in H_1 (L ; \Z/ q\Z)$
 for some $q>1$ instead of $\partial d = 0 \in H_1 (L ; \Z)$, one concludes that the reduction modulo $q$ of  the integer $GW_d^r(X,L ; \underline{x}, \underline{y} , J)$
 neither depends on the choice of $\underline{x}, \underline{y}$ nor on the generic choice of $J$. 
   
\begin{cor}
\label{corlower}
Under the hypothesis of Theorem \ref{theo1disc}, the cardinality of the set \\
${\cal M}_{r,s}^d (X,L;  J , \underline{x} , \underline{y})$ is bounded from below by $\vert GW_d^r (X,L) \vert$. $\square$
\end{cor}

\begin{prop}
\label{propSFT}
Under the hypothesis of Theorem $2.1$, assume that $L$ is a Lagrangian sphere and that $r=1$. Then, the lower bounds given by Corollary \ref{corlower}
are sharp, achieved by any generic almost-complex structure with a very long neck near $L$. Moreover, \\
$(-1)^{[d] \circ [L]} GW_d^1 (X,L) \leq 0$. 

When $L$ is a torus and $r=1$, $GW_d^1 (X,L) = 0$, while $GW_d^r (X,L)$ always vanishes, whatever $r$ is, when $L$ is of genus greater than one. 
However, in both cases, the lower bounds given by Corollary \ref{corlower}
are still sharp, achieved by any generic almost-complex structure with a very long neck near $L$. 
\end{prop}
The notion of almost-complex structure with very long neck has been introduced by Y. Eliashberg, A. Givental and H. Hofer in \cite{EGH}. 
Also, the intersection index $[d] \circ [L] \in \Z/2\Z$ is well defined even though $d \in H_2 (X , L ; \Z)$ is only a relative homology class, since it
does not depend on the choice of a lift of $d$ in $H_2 (X ; \Z)$.\\

{\bf Proof:}

When the genus of $L$ is greater than one, this result is a particular case of Proposition $1.10$ of \cite{WelsSFT} (see also Proposition $4.4$ of \cite{WelsICM}).
When $L$ is a torus, the proof goes along the same lines as Theorem $1.4$ of \cite{WelsSFT} and we do not reproduce it here. The upshot is
that after splitting $(X, \omega)$ near the flat $L$ in the sense of symplectic field theory, for all $J$-holomorphic disc homologous to $d$ and passing
through $\underline{x}, \underline{y}$, the component in the Weinstein neighborhood of $L$ which contains the unique real point $x$ is a once
punctured disc. At the puncture, the disc is asymptotic to a closed Reeb orbit. The boundary of the disc is thus homologous to a closed geodesic
of the flat torus $L$ and thus not homologous to zero in $L$. As a consequence, for an almost-complex structure with very long neck, none of the
$J$-holomorphic discs homologous to $d$ which pass through $\underline{x}, \underline{y}$ have trivial boundary in homology. 

Finally, when $L$ is a sphere, the proof goes along the same lines as Theorem $1.1$ of \cite{WelsSFT} and we do not reproduce it here.
Again, the upshot is that for an almost-complex structure with very long neck and standard near $L$, all the
$J$-holomorphic discs homologous to $d$ which pass through $\underline{x}, \underline{y}$ have boundary close to a geodesic. 
This boundary thus divides $L$ into two components $L^\pm$ which are both homeomorphic to Lagrangian discs. We can then glue
one such Lagrangian disc to our $J$-holomorphic disc $D$ to get an integral two-cycle which lifts $d$ in homology. The intersection index
of this two-cycle with $L$ equals $m(D) + 1 \mod(2)$, where the $+1$ term is the Euler characteristic of the Lagrangian disc. 
We thus deduce that for every such disc $D$, $m(D) = [d] \circ [L] + 1$, the sharpness and the result. $\square$\\

The following Lemma \ref{lemmadivide}
relates this open Gromov-Witten invariant with the real enumerative invariant introduced in \cite{WelsCRAS}, \cite{WelsInvent}. 

\begin{lemma}
\label{lemmadivide}
Let  $(X , \omega , c_X)$ be a closed real symplectic four-manifold which contains a Lagrangian sphere $L$ in its real locus
$\R X = \text{fix} ( c_X)$. Let $d \in H_2 (X ; \Z)$ be such that $c_1 (X) d >0$ and $(c_X)_* d = -d$ and let $1 \leq r \leq c_1 (X) d -1$ be an odd integer. 
Then, $$\chi^d_r (L) = 2^{s-1} \sum_{d' \in H_2 (X, L ; \Z) \; \vert \; d' - (c_X)_* d' = d} (-1)^{d' \circ [\R X \setminus L]} GW_{d'}^r (X,L) \in  \Z,$$ where $s= \frac{1}{2} ( c_1 (X) d -1 - r)$. 
In particular, $2^{s-1}$ divides $\chi^d_r $.
\end{lemma}
Recall that a real symplectic manifold is a symplectic manifold $(X , \omega)$ together with an antisymplectic involution $c_X$,
see \cite{Vit1}, \cite{WelsCRAS}. The invariant $\chi^d_r $ has been introduced in \cite{WelsCRAS}, \cite{WelsInvent}. When the real
locus $\R X$ is not connected, it is understood that  $\chi^d_r (L)$ denotes the part of the invariant obtained by choosing all
the $r$ real points in the component $L$, see \cite{WelsDuke}, \cite{WelsStrong}. Note that if $r$ does not have the same parity as $ c_1 (X) d -1$,
then all invariants vanish by convention so that the formula of Lemma \ref{lemmadivide} holds true. Note finally that even though
$d' \in H_2 (X, L ; \Z) $, the difference $d' - (c_X)_* d' $ is well defined in $H_2 (X ; \Z)$ since it does not depend on the choice of
a lift of $d'$ in $H_2 (X ; \Z)$. Also, the intersection index of $d'$ with the complement $\R X \setminus L$ is well defined. \\

{\bf Proof:}

Let $J \in \R {\cal J}_\omega$ be generic, see \cite{WelsInvent} and $\underline{x}$ (resp.  $\underline{y}$) be a collection of $r$ distinct points
in $L$ (resp. $s$ pairs of complex conjugated points in $X \setminus \R X$). By definition, $\chi^d_r (L) = \sum_{C \in {\cal R}_d (\underline{x}, \underline{y} , J)}
(-1)^{m(C)} \in \Z$, where $ {\cal R}_d (\underline{x}, \underline{y} , J)$ denotes the finite set of real rational $J$-holomorphic curves homologous to $d$
which contain $\underline{x} \cup \underline{y} $. Each such real rational curve is the union of two holomorphic discs with boundary on $L$, exchanged by $c_X$.
Their relative homology class $d'$ thus satisfies $d' - (c_X)_* d' =d$. Moreover, each of these discs contains $\underline{x}$ and one point of every pair
of complex conjugated points $\{ y_i , \overline{y}_i \}$, $1 \leq i \leq s$. There are $2^s$ ways to choose such a point in every pair and we denote by $Y$
this set of $2^s$ $s$-tuples. Now, Schwartz' reflection associates to every $J$-holomorphic disc $u' : \Delta \to X$ with boundary on $L$
a real rational $J$-holomorphic curve $u : \C P^1 \to X$ such that $u \circ conj = conj \circ u$ and its restriction to the upper hemisphere coincides with $u'$. 
We deduce from this Schwartz' reflection a surjective map
$$S : \bigcup_{ \underline{y}' \in Y} \; \bigcup_{d' \in H_2 (X, L ; \Z) \; \vert \; d' - (c_X)_* d' = d}  {\cal M}_{r,s}^{d'} (X,L; J,  \underline{x} , \underline{y}') \to {\cal R}_d (\underline{x}, \underline{y} , J).$$
By definition, for every $[u,J, \underline{z}, \underline{\zeta}] \in {\cal M}_{r,s}^{d'} (X,L; J,  \underline{x} , \underline{y})$, 
$m(u) + d' \circ [\R X \setminus L] = m(S(u)) \mod (2)$ whereas this map $S$ is $2:1$. 
Hence the result. $\square$\\

\begin{rem}
1) In the case of the quadric ellipsoid, the invariant $\chi^d_r (L)$ has been computed in \cite{WelsCRAS1}, \cite{WelsSFT} for small $r$. One can proceed in the same way and express
$ GW_{d}^r (X,L)$ in terms of a similar invariant defined in the cotangent bundle of the two-sphere and of enumerative invariants computed
by Ravi Vakil in the second Hirzebruch surface, see Theorem $3.16$ of \cite{WelsSFT}. When $r$ is small, the open Gromov-Witten invariant
of  the cotangent bundle of the two-sphere is easy to compute, see Lemma $3.5$ of  \cite{WelsSFT}, but for larger $r$, such a computation is not known yet.
Note that an algorithm to compute $\chi^d_r (L)$ for every $r$ has been proposed in  \cite{Shu}. 

2) The last part of Lemma \ref{lemmadivide} actually provides a stronger congruence than the one I already established in \cite{WelsCRAS1}, \cite{WelsSFT} using
symplectic field theory, see Theorem $2.1$ of \cite{WelsSFT} or Theorem $1.4$ of \cite{WelsICM}. 

3) Lemma \ref{lemmadivide}
indicates an obstruction to get similar results for higher genus membranes. Indeed, we know that in a simply connected real projective surface, there is one and only one smooth real 
curve (of genus $\frac{1}{2}(L^2 - c_1 (X)L +2)$) in the linear system of an ample real line bundle $L$, which pass through a real collection of $\frac{1}{2}(L^2 + c_1 (X)L)$ points,
whatever this collection is. Now, if there were an analogous open Gromov-Witten invariant obtained by just counting membranes, we would deduce
for configurations with $s$ pairs of complex conjugated points that $2^{s-1}$ divides one, a contradiction. Already for genus zero membranes with two boundary components
one sees such an obstruction. For instance, through say six real points and one complex point in the quadric ellipsoid, there could be, depending on the position of the points,
either one or zero such genus zero membranes with two boundary components homologous to the hyperplane section in $H_2 (X , L ; \Z)$. The reason is that applying Schwartz' reflection
to such a membrane, one gets a dividing real algebraic curve of bidegree $(2,2)$ in the quadric passing through six real points and two complex conjugated ones. But depending
on the position of these points we know that the unique real curve passing through these points can be either dividing or non-dividing, a contradiction. 
\end{rem}

When the Lagrangian surface $L$ is not orientable, we only obtain the following quite weaker result at the moment. 

\begin{theo}
 \label{theo1discno}
 Let $(X , \omega)$ be a closed connected symplectic four-manifold and $L \subset X$ be a closed Lagrangian surface
 homeomorphic to the real projective plane. Let $d \in H_2 (X , L ; \Z)$ be such that $\mu_L (d) >0$
 and $\partial d \neq 0 \in H_1 (L ; \Z/2\Z)$. Let $r,s \in \N$ be such that $r + 2s = \mu_L (d) - 1$ and 
 $\underline{x}$ (resp.  $\underline{y}$) be a collection of $r$ (resp. $s$) distinct points
of $L$ (resp. $X \setminus L$). Let finally $J \in {\cal J}_\omega$ be generic.
 Then, the reduction modulo $2$ of the integer $GW_d^r(X,L ; \underline{x}, \underline{y} , J)$
 neither depends on the choice of $\underline{x}, \underline{y}$ nor on the generic choice of $J$.
  \end{theo}

This Theorem \ref{theo1discno} thus only provides an invariant $GW_d^r(X,L) \in  \Z/2\Z$. 

\subsection{Proof of Theorems \ref{theo1disc} and Theorem \ref{theo1discno}}
\label{subsecproof1disc}

Let $J_0$ and $J_1$ be two generic elements of ${\cal J}_\omega$. It suffices to prove that $GW_d^r(X,L ; \underline{x}, \underline{y} , J_0)
= GW_d^r(X,L ; \underline{x}, \underline{y} , J_1)$, since the group of symplectic diffeomorphisms of $(X , \omega)$ which preserve $L$
acts transitively on collection of distinct points $(\underline{x}, \underline{y}) \in L^r \times (X \setminus L)^s$. Let $\gamma : t \in [0,1] \mapsto J_t \in {\cal J}_\omega$ 
be a generic path such that $\gamma (0) = J_0$ and $ \gamma (1) = J_1$. Denote by ${\cal M}_\gamma = \pi^{-1} (Im(\gamma))$ the one-dimensional
submanifold of ${\cal M}_{r,s}^d (X,L; \underline{x} , \underline{y})$ and by $\pi_\gamma : {\cal M}_\gamma \to [0,1]$ the associated projection. 
From Propositions \ref{propnondeg} and \ref{propreducible}, $\pi_\gamma$ has finitely many critical points, all non-degenerated, 
which correspond to simple discs with a unique ordinary cusp on their boundary. All the other points of  ${\cal M}_\gamma $ correspond to immersed discs. 

\begin{lemma}
\label{lemmaimm}
 Let $(X , \omega)$ be a closed symplectic four-manifold and $L \subset X$ be a closed Lagrangian surface.
 Let $t \in [0,1] \mapsto J_t \in {\cal J}_\omega$ be a generic path of tame almost-complex structures. Let $u_t : \Delta \to X$, $t \in [0,1]$,
 be a continuous family of $J_t$-holomorphic immersions such that $u_t (\partial  \Delta ) \subset L$. Then, the intersection index
 $[u_t ( \stackrel{\circ}{\Delta})] \circ [L] \in \Z/2\Z$ does not depend on $t \in [0,1]$.
\end{lemma}

{\bf Proof:}

Since these maps are immersions, $L$ is Lagrangian and $ [0,1]$ compact, there exists $\epsilon >0$ such that for every $t \in [0,1]$,
$[u_t ( \stackrel{\circ}{\Delta})] \circ [L] = [u_t ( \Delta (1-\epsilon))] \circ [L]  \in \Z/2\Z$ and $u_t ( \partial \Delta (1-\epsilon))
\cap L = \emptyset$, where $\Delta (1 - \epsilon) = \{ z \in \C \; \vert \; \vert z \vert \leq 1 - \epsilon \}$. The sum of $u_0 ( \Delta (1-\epsilon))$,
the chain $(t, z) \in [0,1] \times \partial \Delta (1-\epsilon) \mapsto (t,  u_t ( z))$ and $- u_1 ( \Delta (1-\epsilon))$ defines an integral two-cycle homologous
to zero. Its intersection index modulo two with $L$ equals $[u_1 ( \stackrel{\circ}{\Delta})] \circ [L] - [u_0 ( \stackrel{\circ}{\Delta})] \circ [L]$,
hence the result. $\square$

\begin{lemma}
\label{lemmacrit}
Let $ [u_{t_0},J_{t_0}, \underline{z}_{t_0}, \underline{\zeta}_{t_0}] \in {\cal M}_\gamma$ be a critical point of $\pi_\gamma$ which is a local
maximum (resp. minimum). There exists a neighborhood $W$ of $ [u_{t_0},J_{t_0}, \underline{z}_{t_0}, \underline{\zeta}_{t_0}]$ in  ${\cal M}_\gamma$
and $\eta >0$ such that for every $t \in ]t_0 - \eta , t_0[$ (resp. $t \in ]t_0 , t_0 +  \eta[$), $\pi_\gamma^{-1} (t) \cap W$ consists of two points
$[u^\pm_{t},J_{t}, \underline{z}^\pm_{t}, \underline{\zeta}^\pm_{t}]$ such that $m(u_t^+) \neq m(u_t^-) \in \Z/2\Z$ and 
for every $t \in ]t_0 , t_0 +  \eta[$ (resp. $t \in ]t_0 - \eta , t_0[$), $\pi_\gamma^{-1} (t) \cap W$ is empty. 
\end{lemma}

{\bf Proof:}

This Lemma \ref{lemmacrit} is the strict analog of Proposition $2.16$ of \cite{WelsInvent} and since its proof also goes
along the same lines, we just recall a sketch of it for the reader's convenience. Let $B$ be a small enough ball centered at the critical
value of $u_{t_0}$ in $X$. The one-parameter family of pseudo-holomorphic discs parameterized by $W$ provides by restriction
a one-parameter family of half-discs in $B$, all immersed except the one corresponding to $t_0$. These half-discs are
parameterized by the set $\{ z \in \Delta \; \vert \;  \Im ( z) \geq 0 \}$. In the space of pseudo-holomorphic half-discs of $B$,
the ones having an ordinary cusp on $L \cap B$ form a codimension one submanifold both sides of which are distinguished
by the intersection index $m(u)$ in $L$. Indeed, representatives of both sides are given by the two perturbations of the
ordinary cusp defined by $t \in \R \mapsto (t^2 , t^3 \pm t)$. The point is then that  the path given by $W$
is transverse to this codimension one wall of non-immersed half-discs, see  \cite{WelsInvent}. $\square$\\

When the one-dimensional cobordism ${\cal M}_\gamma$ is compact, Lemmas \ref{lemmaimm} and \ref{lemmacrit} imply
Theorem \ref{theo1disc}, since as $t$ varies in $[0,1]$, Lemma \ref{lemmaimm} guarantees the invariance of 
$GW_d^r(X,L ; \underline{x}, \underline{y} , J_t)$ between two critical values of $\pi_\gamma$ whereas Lemma \ref{lemmacrit} 
guarantees its invariance while crossing the critical values. 

\begin{figure}[htb]
\begin{center}
\includegraphics[height=4cm]{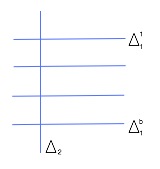}
\end{center}
\caption{Antenna with $b=4$ branches}
\label{figantenna}
\end{figure}

When ${\cal M}_\gamma$ is not compact, its Gromov compactification is given by the following Proposition \ref{propreducible}.
We use here the terminology of Abramovich and Bertram, see Remark 3.2.2 of \cite{AbraBert}, and call an antenna the nodal 
disc obtained by attaching $b$ discs $\Delta_1^1, \dots , \Delta_1^b$ to a single one $\Delta_2$ along $b$ distinct points of 
its boundary $\partial \Delta_2$, see Figure \ref{figantenna}.

\begin{prop}
\label{propreducible}
Let $L$ be a closed orientable Lagrangian surface of a closed connected symplectic four-manifold $(X , \omega)$. Let $d \in H_2 (X , L ; \Z)$ be such that $\mu_L (d) >0$ and
 $r,s \in \N$ such that $r + 2s = \mu_L (d) - 1$. Let $\underline{x}$ (resp.  $\underline{y}$) be a set of $r$ (resp. $s$) distinct points
of $L$ (resp. $X \setminus L$).
Let $t \in [0,1] \mapsto J_t \in {\cal J}_\omega$ be a generic path. Then, every $J_t$-holomorphic element of ${\cal P}_{r,s} (X,L) $,
homologous to $d$ and containing $\underline{x} \cup \underline{y}$, $t \in [0,1] $, is simple. Likewise, every $J_t$-holomorphic
stable disc homologous to $d$ which contains $\underline{x} \cup \underline{y}$, $t \in [0,1] $, is either irreducible, or its image contains
exactly two irreducible components which are both immersed discs transversal to each other. In the latter case, the source disc is an antenna
made of simple discs and the $b$ branches of the antenna are mapped to the same disc, which is of Maslov index zero whenever $b >1$.
\end{prop}

{\bf Proof:}

Let us assume that a $J_t$-holomorphic element of ${\cal P}_{r,s} (X,L) $
homologous to $d$ and containing $\underline{x} \cup \underline{y}$, $t \in [0,1] $, is not simple. From Theorem \ref{theodecomp},
such a $J_t$-holomorphic disc $u : \Delta \to X$ splits into simple discs $u_{\cal D} : \Delta \to X$, so that
$d = \sum_{\cal D} m_{\cal D} (u_{\cal D})_* [\Delta] \in H_2 (X,L ; \Z)$. Let us denote by $({\cal D}_i)_{i \in I}$ the connected
components of  the complement $\Delta \setminus  {\cal G} (u)$ given by Theorem \ref{theodecomp} and set $d_i = (u_{{\cal D}_i})_* [\Delta]$.
From Theorem \ref{theodim}, the $\mu_L (d_i)$ are non-negative.
But the union of these simple discs has to contain the points $\underline{x} \cup \underline{y}$, so that from Theorem \ref{theodim},
$\sum_{i \in I} (\mu_L (d_i) - 1) \geq \# (\underline{x} \cup \underline{y}) -1 \geq \mu_L (d) - 2$.
Since $ \sum_{i \in I} m_{{\cal D}_i} \mu_L (d_i) = \mu_L (d)$ and these Maslov indices are even, we deduce that $\# I =2$.
Hence, either $m_{{\cal D}_1} = m_{{\cal D}_2} = 1$ or $\mu_L (d_i) = 0$  for some $i \in \{ 1, 2 \}$. In both case, one of the discs, say $u_{{\cal D}_2}$, 
together with its incidence conditions, has vanishing Fredholm index, while the other one, $u_{{\cal D}_1}$,
has index $-1$ and a common edge with $u_{{\cal D}_2}$. Perturbing the almost-complex structure on the image of $u_{{\cal D}_2}$,
one observes that the latter condition is of positive codimension while it is independent of the former Fredholm $-1$ condition. As a consequence,
such a non-simple disc cannot appear over a generic path of almost-complex structures $(J_t)_{t \in [0,1] }$. 

Let us assume now that there exists a $J_t$-holomorphic disc $u$ given by Theorem \ref{theocompact} which contains $\underline{x} \cup \underline{y}$
and is homologous to $d$, $t \in ]0,1[$. The preceding arguments again lead to the fact that such a $J_t$-holomorphic disc contains in its
image exactly two irreducible components, both being discs. Indeed, from Theorem \ref{theodim}, the Fredholm index of a simple disc in a class 
$d_1 \in H_2 (X , L ; \Z)$ is 
$\mu_L (d_1) - 1$ whereas
the index of a pseudo-holomorphic sphere in a class $d_2 \in H_2 (X  ; \Z)$ is $\mu_L (d_2) - 2$. As a consequence, the index of
a stable disc with $\alpha$ disc-components  and $\beta$ spherical components in its image equals
$\mu_L (d') - \alpha - 2\beta$, where $d' \in H_2 (X , L ; \Z)$ denotes the total homology class of these image components. We deduce that
$\mu_L (d') - \alpha - 2\beta  \geq \# (\underline{x} \cup \underline{y}) -1$, so that either $\alpha = 2$, $\beta = 0$ and $\mu_L (d') = \mu_L (d)$,
or  $\alpha =0$ and $\beta = 1$. The latter is however excluded since 
$r$ is odd so that  the index of a $J_t$-holomorphic sphere passing through $\underline{x} \cup \underline{y}$ is
then less than $-1$ in this case. Hence, the image of $u$ contains two irreducible components of classes $d_1, d_2 \in H_2 (X  ; \Z)$
such that $\mu_L (d_1) + \mu_L (d_2) = \mu_L (d)$. 

If $\mu_L (d_1)$ and $ \mu_L (d_2) $ are positive, we deduce that the stable disc at the source of $u$ also just contains
two irreducible components and as in the first part of this proof, the restriction of $u$ to both components is simple. These discs are 
immersed from Proposition \ref{propnondeg} and transverse to each other, which follows along the same lines as Proposition $2.11$ of \cite{WelsInvent}. 

If on the contrary one of these Maslov indices, say $\mu_L (d_1)$, vanishes, then $u$ may contain several components at the source which
have the same image homologous to $d_1$. In this case though, the source disc is an antenna, see Figure \ref{figantenna}.
Indeed, denote by $\Delta_2^1$ the unique disc at the source of $u$ whose image is homologous to $d_2$ and by $\Delta_1^1, \dots , \Delta_1^b$
the $b$ stable discs attached to $\Delta_2^1$ whose image have vanishing Maslov index. The neighborhood of $u(\Delta_2^1)$ is standard in $X$,
and the immersion $u \vert_{\Delta_2^1}$ extends to an immersion from a neighborhood $V_2$ of the zero section in the rank two trivial symplectic bundle
equipped with a rank one subbundle over the boundary of $\Delta_2^1$ which has Maslov index $ \mu_L (d_2) $. Likewise, the disc in the image of 
$u$ which has vanishing Maslov index is the image of a $J_t$-holomorphic immersion from a disc $\Delta_1$, which smoothly extends to an immersion 
from a  neighborhood $V_1$ of the zero section in the rank two trivial symplectic bundle equipped with a rank one trivial subbundle over the boundary of $\Delta_1$.
These rank one subbundles over the boundaries of $\Delta_1$ and $\Delta_2^1$ get mapped to $L$. Let us glue $\Delta_1$ and $\Delta_2^1$ along $b$ distinct
points of their boundaries to get a reducible disc $\widetilde{\Delta}$ and glue $V_1$ and $V_2$ along neighborhoods of these points in order to get a manifold $V$
together with an immersion $\tilde{u}_X : V \to X$ and a holomorphic map $\phi : \Delta_1^1 \cup \dots  \cup \Delta_1^b  \cup \Delta_2^1 \to \widetilde{\Delta}$ 
such that $u = \tilde{u}_X \circ \phi$.

Let now $\Sigma_2$ be the second rational ruled surface which contains an exceptional sphere $E$ with self-intersection index $-2$ which is invariant
under the antiholomorphic involution $conj : \Sigma_2 \to \Sigma_2$ whose fixed locus $\R \Sigma_2$ is homeomorphic to a torus. Let $F$ be an irreducible
immersed real rational curve of $\Sigma_2$ which intersects $E$ transversely in $b$ distinct real points and  satisfies $\langle c_1 (\Sigma_2) , F \rangle \geq \mu_L (d_2)$.
Blowing up $\Sigma_2$ along some real points of $F$ if necessary, there exists an immersion $\tilde{u}_{\Sigma_2} : V \to \Sigma_2$ such that 
$\tilde{u}_{\Sigma_2}\vert_{\widetilde{\Delta}} : \widetilde{\Delta} \to E \cup F$ is holomorphic. Suppose then that the $J_t$-holomorphic map 
$u : \Delta_1^1 \cup \dots  \cup \Delta_1^b  \cup \Delta_2^1 \to X$ deforms as $J_{t'}$-holomorphic maps for $t'$ close to $t$, from an irreducible disc $\Delta$ to $X$
in such a way that the contributions of these maps to $GW_d^r(X,L ; \underline{x}, \underline{y} , J_{t'})$ is non-trivial. This implies that the holomorphic map
$\phi : \Delta_1^1 \cup \dots  \cup \Delta_1^b  \cup \Delta_2^1 \to V$ deforms as finitely many $\tilde{u}_X^{-1} (J_{t'})$-holomorphic maps from $\Delta$ to $V$
with boundary on the rank one subbundle and thus that $\tilde{u}_{\Sigma_2} \circ \phi : \Delta_1^1 \cup \dots  \cup \Delta_1^b  \cup \Delta_2^1 \to \Sigma_2$
deforms as finitely many $\tilde{u}_{\Sigma_2} \circ \tilde{u}_X^{-1} (J_{t'})$-holomorphic maps from $\Delta$ to $\Sigma_2$, passing through 
$\tilde{u}_{\Sigma_2} \circ \tilde{u}_X^{-1} (\underline{x} \cup \underline{y})$, have boundary on $\R \Sigma_2$ and non-trivial contribution to 
$GW_d^r(\Sigma_2 ,L ; \tilde{u}_{\Sigma_2} \circ \tilde{u}_X^{-1} (\underline{x}),\tilde{u}_{\Sigma_2} \circ \tilde{u}_X^{-1} ( \underline{y}) , \tilde{u}_{\Sigma_2} \circ \tilde{u}_X^{-1} (J_{t'}))$.
But the almost-complex structure $\tilde{u}_{\Sigma_2} \circ \tilde{u}_X^{-1} (J_{t})$ can be deformed to an almost complex structure $J$ of $\tilde{u}_{\Sigma_2} (V)$
which extends by Schwartz reflection to an almost-complex structure on the neighborhood of $E \cup F$ and such that $\tilde{u}_{\Sigma_2} \circ \phi$
remains holomorphic along this path of almost-complex structures. From Lemma \ref{lemmacrit} follows that the $J$-holomorphic map $\tilde{u}_{\Sigma_2} \circ \phi$ should 
also deforms as nearby holomorphic maps $\Delta \to \Sigma_2$ passing through 
$\tilde{u}_{\Sigma_2} \circ \tilde{u}_X^{-1} (\underline{x} \cup \underline{y})$ and with boundary on $\R \Sigma_2$. From Proposition $3.2.1$
of \cite{AbraBert} follows that $ \Delta_1^1 \cup \dots  \cup \Delta_1^b  \cup \Delta_2^1$ is an antenna and that $u$ restricted to any of its components is simple. 
$\square$\\

 Let $([u_1 ,J_{t_0}, \underline{z}_1, \underline{\zeta}_1]$ , $[u_2 ,J_{t_0}, \underline{z}_2 , \underline{\zeta}_2])$
be a pair of discs given by Proposition \ref{propreducible}, $D_i = Im(u_i)$, and $d_i = (u_i)_* (\Delta) \in H_2 (X , L ; \Z)$, $i \in \{ 1,2 \}$, so that $bd_1 + d_2 = d$.
From the hypothesis we know that $b[\partial d_1] = -[\partial d_2] \in H_1 (L, \Z)$, so that $[\partial d_1] \circ [\partial d_2] =0$. 
Let us equip $L$ with an orientation and denote by $R^+$ (resp. $R^-$) the number of positive (resp. negative) intersection points
between $\partial D_1$ and $\partial D_2$, the latter being canonically oriented by the complex structure. We deduce that
$R^+ = R^- $. Our main observation is then the following, which we formulate in a more general setting since it will be useful
as well in the next paragraph. 

\begin{prop}
\label{propmain}
 Let $(X , \omega)$ be a closed connected symplectic four-manifold and $L \subset X$ be a closed Lagrangian surface.
 Let $([u_1 ,J_{0}, \underline{z}_1, \underline{\zeta}_1] , [u_2 ,J_{0}, \underline{z}_2 , \underline{\zeta}_2]) \in {\cal M}_{r_1,s_1}^{d_1} (X,L; J_0,  \underline{x}_1 , \underline{y}_1) 
 \times {\cal M}_{r_2,s_2}^{d_2} (X,L; J_0,  \underline{x}_2 , \underline{y}_2)$ be a pair of immersed pseudo-holomorphic discs
 transversal to each other. Assume that $r_1 + 2s_1 = \mu_L (d_1)$, $r_2 + 2s_2 = \mu_L (d_2) -1$ and that $L$ is orientable and oriented 
 in the neighborhood of $u_1 (\partial \Delta)$. Let $(J_\lambda)_{\lambda \in ]-\epsilon , \epsilon[}$ be a path transversal to the Fredholm $-1$
 projection $\pi_1 : {\cal M}_{r_1,s_1}^{d_1} (X,L;   \underline{x}_1 , \underline{y}_1) \to {\cal J}_\omega$. Then, as soon as $\epsilon$ is small enough, for every intersection point $w$ of
  $u_1 (\partial \Delta) \cap u_2 (\partial \Delta)$, the pair $([u_1 ,J_{0}, \underline{z}_1, \underline{\zeta}_1] , [u_2 ,J_{0}, \underline{z}_2 , \underline{\zeta}_2])$
  deforms by perturbation of $w$ to exactly one disc in ${\cal M}_{r_1 + r_2,s_1+s_2}^{d_1+d_2} (X,L; J_\lambda,  \underline{x}_1 \cup \underline{x}_2, \underline{y}_1 \cup \underline{y}_2) $
  for positive $\lambda$ and does not deform in ${\cal M}_{r_1 + r_2,s_1+s_2}^{d_1+d_2} (X,L; J_\lambda,  \underline{x}_1 \cup \underline{x}_2, \underline{y}_1 \cup \underline{y}_2) $
  for negative $\lambda$, or vice versa. Moreover, the values of  $\lambda$ for which such a deformation holds only depend on the sign of the local
  intersection index  $u_1 (\partial \Delta) \circ u_2 (\partial \Delta)$ at $w$ in $L$.
\end{prop}

In the statement of Proposition \ref{propmain}, it is understood that the points $\underline{x}_1, \underline{x}_2, \underline{y}_1, \underline{y}_2$ are not
on singular points of the stable disc. Also, the last part of the statement means that the values of $\lambda$ for which the deformation holds 
do not  depend on the Maslov index $ \mu_L (d_2)$, the value of $r_2$ or on the chosen disc $[u_2 ,J_{0}, \underline{z}_2 , \underline{\zeta}_2]$ and intersection point $w$. \\

{\bf Proof:}

Since all nearby discs to the stable one are immersed, we may, without loss of generality,
deform $(J_\lambda)_{\lambda \in ]-\epsilon , \epsilon[}$ or $J_0$ in the image of $\pi_1$ as long as
$([u_1 ,J_{0}, \underline{z}_1, \underline{\zeta}_1] , [u_2 ,J_{0}, \underline{z}_2 , \underline{\zeta}_2])$ remains $J_0$-holomorphic,
as follows from symplectic isotopy.
We may then assume that $J_0$ is standard near $\underline{x}_2, \underline{y}_2$ and blow up these
points to restrict ourselves to the case where $r_2 = s_2 = 0$. 
We are then going to deduce Proposition \ref{propmain} from Proposition $2.14$ of \cite{WelsInvent}. 
Let $b,f \in \R P^1$ and $F' = \{ b \} \times \C P^1$, $B' = \C P^1  \times  \{ f \} $ be the two associated real
rational curves of $ \C P^1 \times  \C P^1$. These curves transversely intersect at one real point. We may
add a finite number of pairs of complex conjugated sections 
$B_i = \C P^1  \times  \{ f_i \} $ 
to $F'$, 
where $f_i \in C P^1 \setminus \R P^1$ and
blow up some pairs of complex conjugated points on $F'$ together with one real point on $B'$ in order to get, after perturbation, 
a real projective surface $(Y, c_Y)$ whose real locus is a once blown-up torus, together with smooth real rational curves $B,F$
intersecting transversely at one real point.  This number of sections and blown-up points can moreover
be chosen such that $c_1 (Y) [F] = \mu_L (d_1)$,
where $d_i = (u_i)_* (\Delta) \in H_2 (X , L ; \Z)$, $i \in \{ 1,2 \}$.
Let $F^+$ be one hemisphere of $F$, so that its interior is a connected component of $F \setminus \R F$.
It induces an orientation on $\partial F^+ =  \R F$. Let us choose an orientation of  $\R Y$ near $\R F$
and denote by $B^+$ (resp. $B^-$) the hemisphere of $B$ such that $[\partial F^+] \circ [\partial B^\pm] = \pm 1$ locally at the unique intersection point. 
Let $U^+$ (resp. $U^-$) be a neighborhood of $F^+ \cup B^+$ (resp. $F^+ \cup B^-$) in $Y$.  
Choose $r_1$ 
distinct points on $\R F \setminus B$ 
and $s_1$ 
distinct pairs of complex conjugated points on $F \setminus \R F$,
where $r_1 + 2s_1 = c_1 (Y) [F] $.
Denote by ${\cal J}_Y$ (resp. $\R {\cal J}_Y$)
the space of almost-complex structures of $Y$ tamed by a chosen symplectic form $\omega_Y$ (resp. for which $c_Y$
is antiholomorphic) and by $J'_0$ the given complex structure of $Y$. Let us denote by ${\cal J}_Y^+$ (resp. $\R {\cal J}_Y^+$)
the codimension one subspace of ${\cal J}_Y$ (resp. $\R {\cal J}_Y$) for which $F^+$ deforms as a pseudo-holomorphic disc
with boundary on $\R Y$ passing through the chosen points on $F^+$. This space contains $J'_0$ and coincides with the
space for which the pair $F^+ \cup B^+$ or $F^+ \cup B^-$ deforms as a pair of pseudo-holomorphic discs
with boundary on $\R Y$ passing through the chosen points. Let $z$ be a smooth point of $F^+ \setminus \partial F^+$ away from
the already chosen ones and $(J'_\lambda)_{\lambda \in ]-\epsilon , \epsilon[}$ be a path of $\R {\cal J}_Y$ which coincides with
$J'_0$ outside of a neighborhood of $z$ and is transversal to ${\cal J}_Y^+$, see \cite{WelsInvent}. 
From Proposition $2.14$ of \cite{WelsInvent}, we may assume that for $\lambda \in ]0 , \epsilon [$ (resp. $\lambda \in ]-\epsilon , 0[$),
$F^+ \cup B^+$ (resp. $F^+ \cup B^-$) deforms as a unique $J'_\lambda$-holomorphic disc with boundary on $\R Y$ passing through the chosen points,
whereas for $\lambda \in ]-\epsilon , 0[$ (resp. $\lambda \in ]0 , \epsilon [$), it has no $J'_\lambda$-holomorphic deformation
with boundary on $\R Y$ passing through the chosen points. These deformations, being close to $F^+ \cup B^\pm$, are all immersed
and thus regular from Proposition \ref{propimmersion}. As a consequence, we deduce by symplectic isotopy that
the latter result neither depends on the choice of the transversal path $(J'_\lambda)_{\lambda \in ]-\epsilon , \epsilon[}$, nor even of the
crossing point $J'_0$  in ${\cal J}_Y$ which makes $F^+ \cup B^\pm$ holomorphic, since this set is connected. Moreover, this result does not depend on
the specific position of the chosen points on $F $, $r_1$ being fixed, since this set is again connected. 

Now, for every positive (resp. negative) intersection point $w_+$ (resp. $w_-$) of $u_1 (\partial \Delta) \cap u_2 (\partial \Delta)$, there exists,
deforming $\omega_Y$ if necessary, a symplectic immersion from $U^+$ (resp. $U^-$) to a neighborhood $V$ of $u_1 ( \Delta) \cup u_2 ( \Delta)$
in $X$ which maps $\R Y \cap U^+$ (resp. $\R Y \cap U^-$) on $L$ and $w$ on $w_+$ (resp. $w_-$). Moreover, there exists holomorphic parameterizations $v_1 :  \Delta \to F^+$
and $v_2 :  \Delta \to B^+$ (resp. $v_2 :  \Delta \to B^-$) such that $u_1 = \phi \circ v_1$ and $u_2 = \phi \circ v_2$. We then fix the points
$v_1 ( \underline{z}_1)$, $v_1 ( \underline{\zeta}_1)$ on $F^+$.
The immersion $\phi$ maps the wall
${\cal J}_Y^+$ on the corresponding codimension one subspace of ${\cal J}_\omega$ made of almost-complex structures for which the
disc $[u_1 ,J_{t_0}, \underline{z}_1, \underline{\zeta}_1]$ deforms as a disc passing through $\underline{x}_1 \cup \underline{y}_1$.
The result now follows from the fact that, up to homotopy, the latter does not depend on the choice of the point $w_+$ (resp. $w_-$), so that
on one side of the codimension one subspace of ${\cal J}_\omega$, the stable disc $([u_1 ,J_{t_0}, \underline{z}_1, \underline{\zeta}_1] , [u_2 ,J_{t_0}, \underline{z}_2 , \underline{\zeta}_2])$
deforms by perturbation of $w^\pm$ to a simple disc while on the other side, it does not, the side only depending on the sign $\pm$. $\square$

\begin{prop}
\label{propmain2}
 Let $(X , \omega)$ be a closed connected symplectic four-manifold and $L \subset X$ be a closed Lagrangian surface.
 Let $([u_1 ,J_{0}] , [u_2 ,J_{0}, \underline{z} , \underline{\zeta}]) \in {\cal M}_{0,0}^{d_1} (X,L; J_0) 
 \times {\cal M}_{r,s}^{d_2} (X,L; J_0,  \underline{x} , \underline{y})$ be a pair of immersed pseudo-holomorphic discs
 transversal to each other. Assume that $ \mu_L (d_1) = 0$, $r + 2s = \mu_L (d_2) -1$ and that $L$ is orientable and oriented 
 in the neighborhood of $u_1 (\partial \Delta)$. Let $(J_\lambda)_{\lambda \in ]-\epsilon , \epsilon[}$ be a path transversal to the Fredholm $-1$
 projection $\pi_1 : {\cal M}_{0,0}^{d_1} (X,L) \to {\cal J}_\omega$. Let $\underline{w}$ be a subset of $b$ points in the intersection 
  $u_1 (\partial \Delta) \cap u_2 (\partial \Delta)$, all of these points having the same local intersection index, so that the contribution
  of  $\underline{w}$ to $u_1 (\partial \Delta) \circ u_2 (\partial \Delta)$ is $\pm b$. Then, as soon as $\epsilon$ is small enough, 
  the pair $([u_1 ,J_{0}] , [u_2 ,J_{0}, \underline{z} , \underline{\zeta}])$
  deforms by perturbation of $\underline{w}$ to exactly one disc in ${\cal M}_{r ,s}^{bd_1+d_2} (X,L; J_\lambda,  \underline{x} , \underline{y}) $
  for positive $\lambda$ and does not deform in ${\cal M}_{r,s}^{bd_1+d_2} (X,L; J_\lambda,  \underline{x}, \underline{y}) $
  for negative $\lambda$, or vice versa. Moreover, the values of  $\lambda$ for which such a deformation holds only depend on the sign of the local
  intersection index  $u_1 (\partial \Delta) \circ u_2 (\partial \Delta)$ at $\underline{w}$ in $L$. If the $b$ points of intersection in $\underline{w}$ 
  do not have the same local intersection index, then the pair $([u_1 ,J_{0}] , [u_2 ,J_{0}, \underline{z} , \underline{\zeta}])$
  does not deform at all in ${\cal M}_{r ,s}^{bd_1+d_2} (X,L; J_\lambda,  \underline{x} , \underline{y}) $ by perturbation of $\underline{w}$, whatever $\lambda$ is.
\end{prop}

{\bf Proof:}

From Proposition \ref{propreducible} we know that if the pair $([u_1 ,J_{0}] , [u_2 ,J_{0}, \underline{z} , \underline{\zeta}])$ is the image of an element
of the Gromov compactification $\overline{\cal M}_\gamma$ of ${\cal M}_\gamma$, then this element is an antenna $u : \Delta_1^1\cup \dots \cup \Delta_1^b\cup \Delta_2^1 \to X$,
where $u$ is an immersion. Given a  subset  $\underline{w}$ of $b$ points in the intersection $u_1 (\partial \Delta) \cap u_2 (\partial \Delta)$, there is up to automorphism a unique such
immersion $u$ which maps the $b$ nodal points of the antenna to $\underline{w}$.
 Let $V$ be a neighborhood of the zero section in the normal bundle of $u$, so that $u$ extends to an immersion $u : V \to X$. 
Let us write $V = V_1^1 \cup \dots \cup V_1^b \cup V_2^1$, where $V_i^j$ is a tubular neighborhood of $ \Delta_i^j$ in $V$. Let us then perturb the pulled back
almost-complex structure $u^{-1}\vert_{V_1^1} (J_0)$ of $V_1^1$ to $u^{-1}\vert_{V_1^1} (J_\lambda)$ where we can assume that the latter does not
depend on $\lambda$ in the intersection $V_1^1 \cap V_2^1$. From Proposition \ref{propmain}, the reducible disc $\Delta_1^1 \cup \Delta_2^1$ deforms to 
an irreducible $u^{-1}\vert_{V_1^1} (J_\lambda)$-holomorphic disc $\Delta_\lambda^1$ for $\lambda >0$ and disappears for $\lambda <0$, or vice-versa. 
Let us perturb then $u^{-1}\vert_{V_1^2} (J_0)$  to $u^{-1}\vert_{V_1^2} (J_\lambda)$ in $V_1^2$ and successively $u^{-1}\vert_{V_1^3} (J_0)$  to $u^{-1}\vert_{V_1^3} (J_\lambda)$,
..., $u^{-1}\vert_{V_1^b} (J_0)$  to $u^{-1}\vert_{V_1^k} (J_\lambda)$, to finally get $u^{-1} (J_\lambda)$ in $V$. We deduce from Proposition \ref{propmain} that
if all the local intersection indices $u (\partial \Delta_1^j) \circ u (\partial \Delta_2^1)$ are the same, the $u^{-1} (J_0)$-holomorphic antenna 
$\Delta_1^1\cup \dots \cup \Delta_1^b \cup \Delta_2^1 $
deforms to a unique and irreducible $u^{-1} (J_\lambda)$-holomorphic disc $\Delta_\lambda$ in $V$ for $\lambda >0$ and disappears for 
$\lambda <0$, or vice-versa. 
Likewise, we deduce from Proposition \ref{propmain} that as soon as two of the local intersection indices 
$u (\partial \Delta_1^j) \circ u (\partial \Delta_2^1)$ differ, 
the $u^{-1} (J_0)$-holomorphic antenna $\Delta_1^1\cup \dots \cup \Delta_1^b \cup \Delta_2^1 $
does not deform to a  $u^{-1} (J_\lambda)$-holomorphic disc $\Delta_\lambda$ in $V$, whatever $\lambda \neq 0$ is. Hence the result. $\square$\\

\begin{rem}
In Remark $2.12$ of \cite{WelsInvent}, I forgot to consider the case where elements of $\R {\cal M}^d (x)$ degenerate to an element of the diagonal
$\Delta$ of Corollary $2.10$. This may happen in the presence of a sphere of vanishing first Chern class as in Proposition \ref{propreducible}.
This case should have been included in Proposition $2.14$ of \cite{WelsInvent} and treated as in Proposition \ref{propmain2}.
\end{rem}

Coming back to the proof of Theorem \ref{theo1disc}, we deduce. 

\begin{cor}
\label{cormain}
Let $([u_1 ,J_{t_0}, \underline{z}_1, \underline{\zeta}_1] , [u_2 ,J_{t_0}, \underline{z}_2 , \underline{\zeta}_2]) \in {\cal M}_{r_1,s_1}^{d_1} (X,L; J_0,  \underline{x}_1 , \underline{y}_1)  \times $ \\
$
{\cal M}_{r_2,s_2}^{d_2} (X,L; J_0,  \underline{x}_2 , \underline{y}_2)$ be a pair of transversal immersed pseudo-holomorphic discs
in the Gromov compactification $\overline{\cal M}_\gamma$ of ${\cal M}_\gamma$ such that $bd_1 + d_2 = d$, $b \geq 1$. Let $R^+$ (resp. $R^-$) be the set of positive (resp. negative) intersection points
between $u_1(\partial \Delta)$ and $u_2 (\partial \Delta)$. There exists a neighborhood $W$ of $([u_1 ,J_{t_0}, \underline{z}_1, \underline{\zeta}_1] , [u_2 ,J_{t_0}, \underline{z}_2 , \underline{\zeta}_2])$ 
in $\overline{\cal M}_\gamma$ and $\eta >0$ such that for every $t \in ]t_0 - \eta , t_0[$ (resp. $t \in ]t_0 , t_0 +  \eta[$), $\pi_\gamma^{-1} (t) \cap W$ consists of
exactly $R^+ \choose b $ (resp. $R^- \choose b $) $J_t$-holomorphic discs, each of which is obtained by perturbation of a different subset of $b$ positive (resp. negative) intersection points
of $u_1(\partial \Delta) \cap u_2 (\partial \Delta)$, or vice-versa. $\square$
\end{cor}

This Corollary \ref{cormain} implies Theorem \ref{theo1disc} in general, since it 
guarantees the invariance of 
$GW_d^r(X,L ; \underline{x}, \underline{y} , J_t)$ while crossing the limit values $\pi_\gamma (\overline{\cal M}_\gamma \setminus {\cal M}_\gamma)$.
Indeed all the elements $[u ,J, \underline{z}, \underline{\zeta}] $ of $W$ given by  Corollary \ref{cormain} have same index $m(u)$
and are thus counted with respect to the same sign by $GW_d^r(X,L ; \underline{x}, \underline{y} , J)$, as soon as $W$ is small enough. $\square$ \\

Finally, Theorem \ref{theo1discno} follows along the same lines as Theorem \ref{theo1disc}. When ${\cal M}_\gamma$ is compact,
 $GW_d^r(X,L)$ is invariant. 
When ${\cal M}_\gamma$ is not compact, its Gromov compactification is given from Proposition \ref{propreducible} by adding to it
pairs of transversal immersed pseudo-holomorphic discs. Let $([u_1 ,J_{t_0}, \underline{z}_1, \underline{\zeta}_1]$ , $[u_2 ,J_{t_0}, \underline{z}_2 , \underline{\zeta}_2])$
be such a pair, $D_i = Im(u_i)$, and $d_i = (u_i)_* (\Delta) \in H_2 (X , L ; \Z)$, $i \in \{ 1,2 \}$.
From the hypothesis we know that either $[\partial d_1]$ or  $[\partial d_2] $ vanishes in $H_1 (L, \Z/2\Z)$.
As a consequence, the discs $([u_1 ,J_{t_0}, \underline{z}_1, \underline{\zeta}_1]$ , $[u_2 ,J_{t_0}, \underline{z}_2 , \underline{\zeta}_2])$
have an even number of intersection points. Now, the analog to Corollary \ref{cormain} tells us that there exists
$R^+, R^- \in \N$, such that $R^+ + R^-$ is this even number of intersection points and for every
$t \in ]t_0 - \eta , t_0[$ (resp. $t \in ]t_0 , t_0 +  \eta[$), $\pi_\gamma^{-1} (t) \cap W$ consists of
exactly $R^+ \choose b$ (resp. $R^- \choose b$) $J_t$-holomorphic discs, each of which is obtained by perturbation of a different 
subset of $b>0$ intersection points
of $u_1(\partial \Delta) \cap u_2 (\partial \Delta)$. Since all these discs are counted with respect to the same sign,
$GW_d^r(X,L ; \underline{x}, \underline{y} , J_t)$ jumps by ${R^+ \choose b} - {R^- \choose b}$ while crossing such a limit value $t_0$.
The result follows from the fact that ${R^+ \choose b} = {R^- \choose b} \mod (2)$, since $R^+ = R^- \mod (2)$. $\square$

\subsection{Higher open Gromov-Witten invariants}

We are now going to extend the results of \S \ref{subsec1disc} by introducing $k$-disc open Gromov-Witten invariants, 
for every positive integer $k$. The advantage is that not only pseudo-holomorphic discs with vanishing boundary will
play a r\^ole here. Let us adopt again the notations of \S \ref{subsec1disc}.
Let $(X , \omega)$ be a closed connected symplectic four-manifold.
Let $L \subset X$ be a closed Lagrangian surface of Maslov class $\mu_L \in H^2 (X,L ; \Z)$ and $k \in \N^*$. 
Let $d \in H_2 (X , L ; \Z)$ be such that $\mu_L (d) \geq k$ and
 $r,s \in \N$ such that $r + 2s = \mu_L (d) - k$. Let finally $\underline{x}$ (resp.  $\underline{y}$) be a set of $r$ (resp. $s$) distinct points
of $L$ (resp. $X \setminus L$). For every $J \in {\cal J}_\omega$ generic, denote by ${\cal M}_{r,s}^{d,k} (X,L; J,  \underline{x} , \underline{y})$
the union $\bigcup_{d_1, \dots , d_k \in H_2 (X , L ; \Z) \; \vert \; d_1+ \dots + d_k = d} {\cal M}_{r,s}^{d_1, \dots , d_k} (X,L; J,  \underline{x} , \underline{y})$,
where ${\cal M}_{r,s}^{d_1, \dots , d_k} (X,L; J,  \underline{x} , \underline{y})$ denotes the finite set of unions of $k$
$J$-holomorphic discs with boundaries on $L$, containing $ \underline{x} , \underline{y}$ and with respective relative homology classes
$d_1, \dots , d_k$. We do not prescribe how these $r+2s$ points $ \underline{x} , \underline{y}$ get distributed among the $k$ discs. 
For every unions of $k$ discs 
$ [u,J, \underline{z}, \underline{\zeta}] = \{ [u_1,J, \underline{z}_1, \underline{\zeta}_1] , \dots , [u_k,J, \underline{z}_k, \underline{\zeta}_k] \}
\in {\cal M}_{r,s}^{d,k} (X,L;  \underline{x} , \underline{y})$, we set $m(u) = \sum_{i=1}^k m (u_i) \in \Z/2\Z$. Then, we set 
$$GW_{d,k}^r(X,L ; \underline{x}, \underline{y} , J) = \sum_{ [u,J, \underline{z}, \underline{\zeta}] \in {\cal M}_{r,s}^{d,k} (X,L; J,  \underline{x} , \underline{y})} (-1)^{m(u)} \in \Z.$$
 \begin{theo}
 \label{theokdiscs}
 Let $(X , \omega)$ be a closed connected symplectic four-manifold, $L \subset X$ be a closed Lagrangian surface
 which we assume to be orientable and $k \in \N^*$. Let $d \in H_2 (X , L ; \Z)$ be such that $\mu_L (d) \geq k$
 and $\partial d = 0 \in H_1 (L ; \Z)$. Let $r,s \in \N$ be such that $r + 2s = \mu_L (d) - k$ and 
 $\underline{x}$ (resp.  $\underline{y}$) be a collection of $r$ (resp. $s$) distinct points
of $L$ (resp. $X \setminus L$). Let finally $J \in {\cal J}_\omega$ be generic.
 Then, the integer $GW_{d,k}^r(X,L ; \underline{x}, \underline{y} , J)$
 neither depends on the choice of $\underline{x}, \underline{y}$ nor on the generic choice of $J$. 
 The same holds true modulo two when $L$ is homeomorphic to the real projective plane and $\partial d \neq 0 \in H_1 (L ; \Z/2\Z)$. 
 \end{theo}
 The invariant provided by Theorem \ref{theokdiscs} can be denoted without ambiguity by $GW_{d,k}^r (X,L) \in  \Z$,
 it is also left invariant under deformation of the symplectic form $\omega$,
 whereas it depends in general of $d$ and $r$. Examples where this invariant can be non-trivial are given by toric fibers of toric surfaces,
 already when counting the finite union of Maslov two discs for $r=k$, see Proposition $6.1.4$ of \cite{BiranCor}. 
\begin{cor}
\label{corlowerk}
Under the hypothesis of Theorem \ref{theokdiscs}, the cardinality of the set \\
${\cal M}_{r,s}^{d,k} (X,L;  J , \underline{x} , \underline{y})$
is bounded from below by $\vert GW_{d,k}^r (X,L) \vert$. $\square$
\end{cor}

{\bf Proof of Theorem \ref{theokdiscs}:}

The proof of Theorem \ref{theokdiscs} goes along the same lines as the one of Theorems \ref{theo1disc} and \ref{theo1discno}.
Let $J_0$ and $J_1$ be two generic elements of ${\cal J}_\omega$, it suffices to prove that $GW_{d,k}^r(X,L ; \underline{x}, \underline{y} , J_0)
= GW_{d,k}^r(X,L ; \underline{x}, \underline{y} , J_1) $. Let $\gamma : t \in [0,1] \mapsto J_t \in {\cal J}_\omega$ 
be a generic path such that $\gamma (0) = J_0$ and $ \gamma (1) = J_1$, ${\cal M}_\gamma = \cup_{ t \in [0,1] } {\cal M}_{r,s}^{d,k} (X,L; J_t,  \underline{x} , \underline{y})$ the
corresponding one-dimensional manifold and $\pi_\gamma : {\cal M}_\gamma \to [0,1]$ the associated projection. 
When ${\cal M}_\gamma$ is compact, Theorem \ref{theokdiscs} again follows from Lemmas \ref{lemmaimm} and \ref{lemmacrit}.
When ${\cal M}_\gamma$ is not compact, its Gromov compactification is given from Proposition \ref{propreducible} by adding to it
$(k+1)$-tuple of transversal immersed pseudo-holomorphic discs. Let $([u_1 ,J_{t_0}, \underline{z}_1, \underline{\zeta}_1]$ , \dots , $[u_{k+1} ,J_{t_0}, \underline{z}_{k+1} , \underline{\zeta}_{k+1}])$
be such a $(k+1)$-tuple, $D_i = Im(u_i)$, and $d_i = (u_i)_* (\Delta) \in H_2 (X , L ; \Z)$, $i \in \{ 1,2 \}$.
One of these $k+1$ discs, say $D_1$, is, together with its incidence conditions, of Fredholm index $-1$ while the other are of 
vanishing Fredholm index. Again, from the hypothesis we deduce that $[\partial d_1] \circ( [\partial d_2] + \dots +  [\partial d_{k+1}]) =0 $,
so that by  equipping $L$ with an orientation and denoting by $R^+$ (resp. $R^-$) the set of positive (resp. negative) intersection points
between $\partial D_1$ and $\partial D_2 \cup \dots \cup \partial D_{k+1}$, we have $R^+ = R^- $. The latter holds true modulo $2$
when $L$ is a real projective plane. Theorem \ref{theokdiscs} then follows from Propositions \ref{propmain} and \ref{propmain2}. $\square$\\

When $L$ is a Lagrangian sphere, these $k$-discs open Gromov-Witten invariants are deduced from the one-disc
invariants introduced in \S  \ref{subsec1disc} in the following way. 

\begin{lemma}
 Let $(X , \omega)$ be a closed connected symplectic four-manifold, $L \subset X$ be a Lagrangian sphere 
 and $k \in \N^*$. Let $d \in H_2 (X , L ; \Z)$ be such that $\mu_L (d) \geq k$
 and $r,s \in \N$ such that $r + 2s = \mu_L (d) - k$. Then, 
 $$GW_{d,k}^r (X,L) = \frac{1}{k!} \sum_{d_1, \dots , d_k  \; \vert \;  \sum d_i = d} \; \sum_{r_1, \dots , r_k  \; \vert \; \sum  r_i = r} {r  \choose r_1, \dots , r_k}  {s  \choose s_1, \dots , s_k}
 \prod_{i=1}^k GW_{d_i}^{r_i}(X,L),  $$
where for $1 \leq i \leq k$, $s_i = \frac{1}{2} (\mu_L (d_i) - r_i - 1)$ and ${r  \choose r_1, \dots , r_k}= {r  \choose r_1}{r-r_1  \choose r_2}\dots{r_{k-1} + r_k  \choose r_{k-1}}$.
\end{lemma}

Note that if $r_i = \mu_L (d_i) \mod (2)$, $GW_{d_i}^{r_i}(X,L) = 0$.\\

{\bf Proof:}

Let $\underline{x}$ (resp.  $\underline{y}$) be a collection of $r$ (resp. $s$) distinct points
of $L$ (resp. $X \setminus L$) and $J \in {\cal J}_\omega$ be generic. Then, $k! \, GW_{d,k}^r (X,L) $
counts the number of ordered unions of $k$ discs $ ([u_1,J, \underline{z}_1, \underline{\zeta}_1] , \dots , [u_k,J, \underline{z}_k, \underline{\zeta}_k]) $
containing $\underline{x} \cup \underline{y}$ and of total homology class $d$. The images of $\underline{z}_i, \underline{\zeta}_i$ under $u_i$
induce a partition of  $\underline{x}$,  $\underline{y}$ into subsets $\underline{x}_i$,  $\underline{y}_i$ of respective cardinalities $r_i$, $s_i$, $1 \leq i \leq k$.
The number of such partitions equals respectively ${r  \choose r_1, \dots , r_k}$ and ${s  \choose s_1, \dots , s_k}$, while 
$ ([u_1,J, \underline{z}_1, \underline{\zeta}_1] , \dots , [u_k,J, \underline{z}_k, \underline{\zeta}_k]) $
provides an element counted by the product $ \prod_{i=1}^k GW_{d_i}^{r_i}(X,L)$. Conversely, for every such partition,
any $k$-tuple of discs counted by the invariant $ \prod_{i=1}^k GW_{d_i}^{r_i}(X,L)$ provides an element counted by $k! \, GW_{d,k}^r (X,L) $ with 
respect to the same sign. Hence the result. $\square$\\

When $L$ is fixed by an antisymplectic involution, the invariant $\chi^d_r$ introduced in \cite{WelsCRAS}, \cite{WelsInvent} counts real rational
curves which thus consist of pairs of complex conjugated discs, the sum of the boundaries of which vanish in homology. However, the invariant
$GW_{d,2}^r(X,L)$ seems to count a much larger number of pairs so that I do not see a relation between the two in general.

\addcontentsline{toc}{part}{\hspace*{\indentation}Bibliography}

\bibliographystyle{abbrv}

\begin{thebibliography}{}

\end{thebibliography}


\begin{thebibliography}{10}

\bibitem{AbraBert}
D. Abramovich and A. Bertram.
\newblock The formula 12=10+2?1 and its generalizations: counting
 rational curves on F2.
 \newblock {\em Advances in algebraic geometry motivated by physics (Lowell, MA, 2000)}, 
 83--88, Contemp. Math., 276, Amer. Math. Soc., Providence, RI,  2001. 

\bibitem{BiranCor}
P.~Biran and O.~Cornea.
\newblock Quantum structures for lagrangian submanifolds.
\newblock {\em Preprint arXiv:0708.4221}, 2007.

\bibitem{BruPui}
E.~Brugall\'e and N.~Puignau.
\newblock Behavior of Welschinger Invariants under Morse Simplifications.
\newblock {\em To appear in Rendiconti del Seminario Matematico della Universitˆ di Padova. Preprint 	arXiv:1203.2773}, 2012.

\bibitem{Cho}
C.-H. Cho.
\newblock Counting real {$J$}-holomorphic discs and spheres in dimension four
  and six.
\newblock {\em J. Korean Math. Soc.}, 45(5):1427--1442, 2008.

\bibitem{EGH}
Y.~Eliashberg, A.~Givental, and H.~Hofer.
\newblock Introduction to symplectic field theory.
\newblock {\em Geom. Funct. Anal.}, (Special Volume, Part II):560--673, 2000.
\newblock GAFA 2000 (Tel Aviv, 1999).

\bibitem{Frauen}
U.~Frauenfelder.
\newblock Gromov convergence of pseudoholomorphic disks.
\newblock {\em J. Fixed Point Theory Appl.}, 3(2):215--271, 2008.

\bibitem{Fuk}
K.~Fukaya.
\newblock Counting pseudo-holomorphic discs in Calabi-Yau 3-folds.
\newblock {\em Tohoku Math. J.}, 63(4):697Ð727,  2011.

\bibitem{FOOO}
K.~Fukaya, Y.-G. Oh, H.~Ohta, and K.~Ono.
\newblock {\em Lagrangian intersection {F}loer theory: anomaly and obstruction.
  {P}art {I}, {P}art {II}}, volume~46 of {\em AMS/IP Studies in Advanced
  Mathematics}.
\newblock American Mathematical Society, Providence, RI, 2009.

\bibitem{Gromov}
M.~Gromov.
\newblock Pseudoholomorphic curves in symplectic manifolds.
\newblock {\em Invent. Math.}, 82(2):307--347, 1985.

\bibitem{HLS}
H.~Hofer, V.~Lizan, and J.-C. Sikorav.
\newblock On genericity for holomorphic curves in four-dimensional
  almost-complex manifolds.
\newblock {\em J. Geom. Anal.}, 7(1):149--159, 1997.

\bibitem{IvashShev}
S.~Ivashkovich and V.~Shevchishin.
\newblock Structure of the moduli space in a neighborhood of a cusp-curve and
  meromorphic hulls.
\newblock {\em Invent. Math.}, 136(3):571--602, 1999.

\bibitem{KatzLiu}
S.~Katz and C.-C.~M. Liu.
\newblock Enumerative geometry of stable maps with {L}agrangian boundary
  conditions and multiple covers of the disc.
\newblock {\em Adv. Theor. Math. Phys.}, 5(1):1--49, 2001.

\bibitem{KwonOh}
D.~Kwon and Y.-G. Oh.
\newblock Structure of the image of (pseudo)-holomorphic discs with totally
  real boundary condition.
\newblock {\em Comm. Anal. Geom.}, 8(1):31--82, 2000.
\newblock Appendix 1 by Jean-Pierre Rosay.

\bibitem{Iac}
V. Iacovino.
\newblock Open Gromov-Witten theory on Calabi-Yau three-folds I.   
\newblock {\em Preprint math.SG/0907.5225}, 2009.

\bibitem{Laz}
L.~Lazzarini.
\newblock Decomposition of a {J}-holomorphic curve.
\newblock {\em Preprint http://www.math.jussieu.fr/~lazzarin/articles.html.}

\bibitem{Laz1}
L.~Lazzarini.
\newblock Existence of a somewhere injective pseudo-holomorphic disc.
\newblock {\em Geom. Funct. Anal.}, 10(4):829--862, 2000.

\bibitem{Liu}
C.-C.~M. Liu.
\newblock Moduli of {J}-holomorphic curves with {L}agrangian boundary
  conditions and open {G}romov-{W}itten invariants for an ${S}^1$-equivariant
  pair.
\newblock {\em Thesis ({P}h.{D}.)-{H}arvard {U}niversity,
  arXiv:math.SG/0210257}, 2002.

\bibitem{McDSal}
D.~McDuff and D.~Salamon.
\newblock {\em {$J$}-holomorphic curves and symplectic topology}, volume~52 of
  {\em American Mathematical Society Colloquium Publications}.
\newblock American Mathematical Society, Providence, RI, 2004.

\bibitem{Shu}
E.~Shustin.
\newblock Welschinger invariants of toric del {P}ezzo surfaces with nonstandard
  real structures.
\newblock {\em Tr. Mat. Inst. Steklova}, 258(Anal. i Osob. Ch. 1):227--255,
  2007.

\bibitem{Smale}
S.~Smale.
\newblock An infinite dimensional version of {S}ard's theorem.
\newblock {\em Amer. J. Math.}, 87:861--866, 1965.

\bibitem{Solomon}
J.~P. Solomon.
\newblock Intersection theory on the moduli space of holomorphic curves with
  {L}agrangian boundary conditions.
\newblock {\em Pr\'epublication math.SG/0606429}, 2006.

\bibitem{Vit1}
C.~Viterbo.
\newblock Symplectic real algebraic geometry.
\newblock {\em Preprint}, 1999.

\bibitem{WelsCRAS}
J.-Y. Welschinger.
\newblock Invariants of real rational symplectic 4-manifolds and lower bounds
  in real enumerative geometry.
\newblock {\em C. R. Math. Acad. Sci. Paris}, 336(4):341--344, 2003.

\bibitem{WelsStrong}
J.-Y. Welschinger.
\newblock Enumerative invariants of strongly semipositive real symplectic
  six-manifolds.
\newblock {\em Pr\'epublication math.AG/0509121}, 2005.

\bibitem{WelsInvent}
J.-Y. Welschinger.
\newblock Invariants of real symplectic 4-manifolds and lower bounds in real
  enumerative geometry.
\newblock {\em Invent. Math.}, 162(1):195--234, 2005.

\bibitem{WelsDuke}
J.-Y. Welschinger.
\newblock Spinor states of real rational curves in real algebraic convex
  3-manifolds and enumerative invariants.
\newblock {\em Duke Math. J.}, 127(1):89--121, 2005.

\bibitem{WelsCRAS1}
J.-Y. Welschinger.
\newblock Invariant count of holomorphic disks in the cotangent bundles of the
  two-sphere and real projective plane.
\newblock {\em C. R. Math. Acad. Sci. Paris}, 344(5):313--316, 2007.

\bibitem{WelsSFT}
J.-Y. Welschinger.
\newblock Optimalit\'e, congruences et calculs d'invariants des vari\'et\'es
  symplectiques r\'eelles de dimension quatre.
\newblock {\em Preprint math.SG/0707.4317}, 2007.

\bibitem{WelsICM}
J.-Y. Welschinger.
\newblock Invariants entiers en g\'eom\'etrie \'enum\'erative r\'eelle.
\newblock In {\em Proceedings of the {I}nternational {C}ongress of
  {M}athematicians, {V}ol. {II} ({H}yderabad, 2010)}, pages 652--678, 2010.
  
  \bibitem{Welsopen6}
J.-Y. Welschinger.
\newblock Open Gromov-Witten invariants in dimension six.
\newblock {\em Math. Ann., published online DOI 10.1007/s00208-012-0883-0}, 2012.

\end{thebibliography}

\vspace{0.7cm}
\noindent 
\textsc{Universit\'e de Lyon \\
CNRS UMR 5208 \\
Universit\'e Lyon 1 \\
Institut Camille Jordan} \\
43 blvd. du 11 novembre 1918 \\
F-69622 Villeurbanne cedex\\
France

\end{document}